\documentclass[conference]{IEEEtran}

\IEEEoverridecommandlockouts

\usepackage{algorithm}
\usepackage{algorithmic}
\usepackage[cmex10]{amsmath}
\usepackage{amsfonts, amssymb}
\usepackage{bm} 
\usepackage{enumitem}
\usepackage[caption=false,font=footnotesize]{subfig}
\usepackage[pdftex]{graphicx}
\usepackage{mathrsfs}
\usepackage[dvipsnames]{xcolor}
\usepackage{mathtools}
\usepackage{float}
\usepackage[noadjust]{cite}

\graphicspath{{./figs/}}
\DeclareGraphicsExtensions{.png,.jpg,.pdf,.eps}

\newcommand{\norm}[1]{\left\lVert#1\right\rVert} 
\newcommand{\real}{\mathbb{R}} 
\newcommand{\ones}{\mathbf{1}}

\newcommand{\trans}{^\top}

\newcommand{\diag}[1]{\mathrm{diag}(#1)}
\newcommand{\rank}[1]{\mathrm{rank}(#1)}

\DeclareMathOperator*{\minimize}{minimize}

\begin{document}

\title{A Data-driven Voltage Control Framework for \\ Power Distribution Systems}

\author{Hanchen Xu, Alejandro~D.~Dom\'{i}nguez-Garc\'{i}a, and Peter~W.~Sauer
\thanks{The authors are with the Department of Electrical and Computer Engineering at the University of Illinois at Urbana-Champaign, Urbana, IL 61801, USA. Email: \{hxu45, aledan, psauer\}@illinois.edu.}
}

\maketitle

\begin{abstract}
In this paper, we address the problem of coordinating a set of distributed energy resources (DERs) to regulate voltage in power distribution systems to desired levels.
To this end, we formulate the voltage control problem as an optimization problem, the objective of which is to determine the optimal DER power injections that minimize the voltage deviations from desirable voltage levels subject to a set of constraints.
The nonlinear relationship between the voltage magnitudes and the nodal power injections is approximated by a linear model, the parameters of which can be estimated in real-time efficiently using measurements.
In particular, the parameter estimation requires much fewer data by exploiting the structural characteristics of the power distribution system.
As such, the voltage control framework is intrinsically adaptive to parameter changes.
Numerical studies on the IEEE 37-bus power distribution test feeder validated the effectiveness of the propose framework.

\end{abstract}

\begin{IEEEkeywords}
distributed energy resource, voltage control, parameter estimation, data-driven, power distribution system.
\end{IEEEkeywords}

\section{Introduction}

Distributed energy resources (DERs), which typically have fast-responding characteristics, could potentially be utilized to regulate voltage in power distribution systems (see, e.g., \cite{VC1, VC2} and references therein). 
Many research efforts have been put in developing coordination schemes for DERs to provide voltage control under the assumption that a model of the system is available.
For example, in \cite{VC1}, the authors proposed a two-stage distributed architecture for voltage control in distribution systems, where the required reactive power injections are determined by each local controller in the first stage, and any deficiency is compensated in the second stage by other nodes providing more reactive power so as to uniformly raise the voltage profiles across the network. 
While the voltage control problem in distribution networks can be approached from a control perspective as in \cite{VC1}, it can also be addressed from an optimization perspective \cite{VC2}.
These approaches could effectively regulate the voltages given an accurate model; however, they may fail in the absence of an accurate model, which is typically the case in practice. 
As such, practical voltage control schemes adaptive to system changes and robust against model errors are still to be investigated.

The ongoing deployment of advanced meters and communications makes it possible to access more detailed operating states of a power distribution system. 
Using measurements of system operating states, data-driven voltage control schemes that utilize sensitivities can be developed as an alternative approach for voltage control.
Existing works have applied data-driven approaches to the estimation of sensitivities, including injection shifting factors, and power transfer distribution factors \cite{MB1, MB2, MB3, MB4, ZJB}. 
The estimated sensitivities haven then been used in contingency analysis, generator re-dispatch, and market clearing \cite{MB1, MB2}. 
The sensitivities estimated from measurement enjoy several nice properties, including adaptivity to system changes and robustness to modeling errors. 

Notwithstanding the abundance of works on voltage control and power system sensitivity estimation, to the best of our knowledge, there is no previous work combining these two approaches with the intent of developing a voltage control framework that is adaptive to system changes and robust to model errors.
In this regard, we propose a data-driven voltage control framework that consists of two fundamental building blocks: (i) an estimator of the system model based on measurements; and (ii) a voltage controller based on the estimated model. 
Specifically, we will take advantage of an approximate linear model---the so-called LinDistFlow model (see, e.g., \cite{33bus, DistFlow})---to simlify the nonlinear relationships between voltage magnitudes and power injections. 
The coefficients of this approximate linear model are essentially the sensitivities of the squared voltage magnitudes with respect to active and reactive power injections. 
These coefficients will be estimated online and as such will be adaptive to reflect system changes.
When the system topology is known, by taking advantage of the structural characteristics of the power distribution system, the estimation of the model parameters requires much fewer data, as compared to conventional sensitivity estimation approaches; this enhances the adaptivity of the estimator.
Using the estimated model, we formulate the voltage control problem as a convex optimization problem that can be solved efficiently to determine the optimal set-points for the DER power injections.
Consequently, the voltage controller that builds upon the estimated model will also be adaptive. 


\section{Preliminaries}
In this section, we provide the power distribution model adopted in this work. Subsequently, we describe the problem of regulating voltage in a distribution network. 

\subsection{Branch Flow Model} \label{sec:bfm}

Consider a radial power distribution network represented by a directed graph  $\mathcal{G} = (\tilde{\mathcal{N}}, \mathcal{E})$, where $\tilde{\mathcal{N}} = \{0, 1, \cdots, N\}$ is the set of buses (nodes), and $\mathcal{E} = \tilde{\mathcal{N}} \times \tilde{\mathcal{N}}$ is set of transmission lines (edges).
Let $\mathcal{L} = \{1, 2, \cdots, L\}$ be the index set of transmission lines. 
Each electrical line $\ell \in \mathcal{L}$ is associated with $e=(i, j) \in \mathcal{E}$, i.e., $i$ is the sending end and $j$ is the receiving end of line $\ell$, with the direction from $i$ to $j$ defined to be positive. 
Let $r_\ell$ and $x_\ell$ denote the resistance and reactance of line $\ell \in \mathcal{L}$.
Define $\bm{r}=[r_1, \cdots, r_L]\trans$ and $\bm{x}=[x_1, \cdots, x_L]\trans$.
Throughout the rest of the paper, we make the same assumption that the ``R-to-X" ratios of line $\ell$ equals to a known constant $\alpha_\ell$, i.e., $r_\ell/x_\ell = \alpha_\ell$, and define $\bm{\alpha} = [\alpha_1, \cdots, \alpha_L]\trans$.
When the network is radial, $\mathcal{G}$ has a tree-like structure. Without loss of generality, assume there is one DER and one load at each bus except bus~$0$, which we assume corresponds to the substation bus, and define $ \mathcal{N} := \tilde{\mathcal{N}} \setminus \{ 0\}=\{1, \cdots, N\}$.
We refer to the DER and the load at bus $i$ as DER $i$ and load $i$, respectively.

Let $\tilde{\bm{M}}=[M_{i \ell}] \in \real^{(N+1) \times L}$ denote the node-to-edge incidence matrix of $\mathcal{G}$, with $M_{i \ell} = 1$ and $M_{j \ell} = -1$ if line $\ell$ starts from bus $i$ and ends at bus $j$, and all other entries equal to zero. 
Let $\bm{M}$ denote the $(N\times L)$-dimensional matrix that results from removing the first row in $\tilde{\bm{M}}$.
Let $\mathcal{P}_i \subseteq \mathcal{L}$ denote the set of lines that form the unique path from bus 0 to bus $i$.
For a radial distribution system, $L=N$, and $\bm{M}$ is invertible.

Let $V_i$ denote the magnitude of the voltage at bus  $i \in \tilde{\mathcal{N}}$, by $V_i$,  and define $\bm{V} = [V_1, \cdots, V_N]\trans$; in the remainder, we assume that $V_0$ is a constant.
In addition, let $p_i^g$ and $q_i^g$ respectively denote the active and reactive power injected by DER~$i$, and define $\bm{p}^g = [p_1^g, \cdots, p_N^g]\trans$, and $\bm{q}^g = [q_1^g, \cdots, q_N^g]\trans$. 
Similarly, let $p_i^d$ and $q_i^d$ respectively denote the active and reactive power demanded by load~$i$, and define $\bm{p}^d = [p_1^d, \cdots, p_N^d]\trans$, and $\bm{q}^d = [q_1^d, \cdots, q_N^d]\trans$. 
Let $\underline{p}_i^g$ and $\overline{p}_i^g$ respectively denote the minimum and  maximum active  power that can be provided by DER~$i$, and define $\bm{\underline{p}}^g = [\underline{p}_1^g, \cdots, \underline{p}_N^g]\trans$ and $\bm{\overline{p}}^g = [\overline{p}_1^g, \cdots, \overline{p}_N^g]\trans$. 
Similarly, let $\underline{q}_i^g$ and $\overline{q}_i^g$ respectively denote the minimum and maximum reactive power that can be provided by DER~$i$, and define $\bm{\underline{q}}^g = [\underline{q}_1^g, \cdots, \underline{q}_N^g]\trans$ and $\bm{\overline{q}}^g = [\overline{q}_1^g, \cdots, \overline{q}_N^g]\trans$.
Let $p_i=p_i^g-p_i^d$ and $q_i=q_i^g-q_i^d$, and define  $\bm{p}=[p_1, \cdots, p_N]\trans = \bm{p}^g-\bm{p}^d$, and $\bm{q}=[q_1, \cdots, q_N]\trans = \bm{q}^g-\bm{q}^d$. 
Let $f_\ell$ denote the active power that flows from the sending end to the receiving end of line $\ell \in \mathcal{L}$, and define $\bm{f} = [f_1, \cdots, f_L]\trans$. 
Let  $\overline{f}_\ell$ denote maximum power flows on line $\ell \in \mathcal{L}$, and define $\overline{\bm{f}} = [\overline{f}_1, \cdots, \overline{f}_L]\trans$.

Let $u_i = V_i^2$, and define $\bm{u} = [u_1, \cdots, u_N]\trans$, $\tilde{\bm{u}} = \bm{u} - u_0 \ones_N$.
Assume that the power distribution system is lossless, then $\bm{f}$ can be directly computed from $\bm{p}$ as follows:
\begin{equation} \label{eq:line_flow}
	\bm{f} = \bm{M}^{-1} \bm{p}.
\end{equation}
Then, the relation between $\bm{u}$, $\bm{p}$, and $\bm{q}$, can be captured by the so-called LinDisfFlow model as follows \cite{DistFlow}:
\begin{equation} \label{eq:LinDistFlow}
    \tilde{\bm{u}} = \bm{R} \bm{p} + \bm{X} \bm{q}, 
\end{equation}
where $\ones_N$ is the $N$-dimensional all-ones vector, and
\begin{equation} \label{eq:RX}
\begin{array}{c}
	\bm{R} = 2 (\bm{M}^{-1})\trans \diag{\bm{r}} \bm{M}^{-1},\\[5pt]
	\bm{X} = 2 (\bm{M}^{-1})\trans \diag{\bm{x}} \bm{M}^{-1},
\end{array}
\end{equation}
where $\diag{\cdot}$ returns a diagonal matrix with the entries of the argument on its diagonals; we refer to the matrices $\bm{R}$ and $\bm{X}$ as the voltage sensitivity matrices.

\subsection{Voltage Control Problem} \label{sec:vcp}
The objective here is to maintain the voltage magnitude at each bus $i$, $i \in \mathcal{N}$, of the distribution system   within a pre-specified interval denoted by $[\underline{V}_i, \overline{V}_i]$.
While a number of means, such as tap change under load transformers and capacitor banks can be utilized to achieve  the aforementioned objective, it is also possible to utilize the DERs present in the distribution system. In this paper, we focus solely in this later mechanism for achieving voltage control. 
Then, the problem is to determine the DER active and/or reactive power injections, $\bm{p}^g$ and $\bm{q}^g$, so that
\begin{itemize} 
\item[\textbf{[C1.]}] the active and reactive power injections from each DER $i$, $i \in \mathcal{N}$, do not exceed its corresponding capacity limits, i.e., $\bm{\underline{p}}^g \leq~\bm{p}^g \leq \bm{\overline{p}}^g$, and $\bm{\underline{q}}^g \leq \bm{q}^g \leq \bm{\overline{q}}^g$;  
\item[\textbf{[C2.]}] the voltage magnitude at each bus $i$, $i \in \mathcal{N}$, is within the pre-specified interval, i.e., $\underline{V}_i \leq V_i \leq \overline{V}_i$; and
\item[\textbf{[C3.]}] the power flow on each line $\ell$, $\ell \in \mathcal{L}$, does not exceed its maximum capacity, i.e., $-\overline{\bm{f}} \leq \bm{f} \leq \overline{\bm{f}}$.
\end{itemize}
Note that while constraint \textbf{C1} is a hard constraint that cannot be violated, constraints \textbf{C2} and \textbf{C3} may be allowed to be violated slightly.
In addition, among all feasible values of $\bm{p}^g$ and $\bm{q}^g$, we would like to select the one that minimizes some cost function, e.g., one that reflects the cost of voltage deviations as well as the cost of active/reactive power provision.

\section{Voltage Control Framework} \label{sec:framework}

In this section, we propose an adaptive data-driven framework for voltage control using DERs. 
We first give an overview on the framework and then present the details of the fundamental building blocks that compromise the framework.

\subsection{Framework Overview}

\begin{figure}[!t]
\centering
\includegraphics[width=3in]{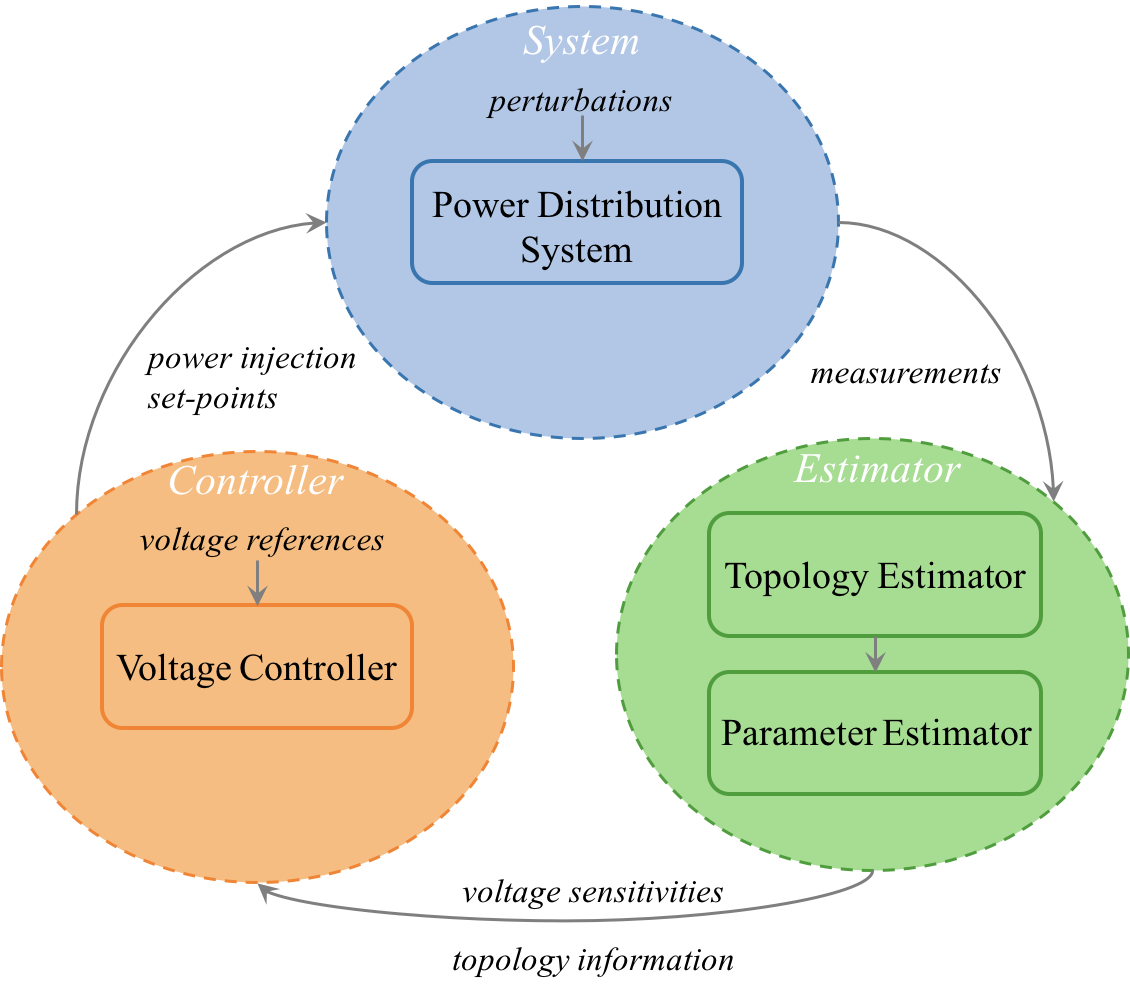}
\caption{Data-driven voltage control framework.}
\label{fig:framework}
\end{figure}

The proposed voltage control framework consists of two components, an estimator and a controller.
The interaction between the different components is illustrated via the block diagram in Fig. \ref{fig:framework}.
The estimator component contains a topology estimator that estimates the topology of the power distribution system (essentially, $\bm{M}$), and a parameter estimator that estimates the line parameters ($\bm{r}$ and $\bm{x}$), using measurements of power injections and voltage magnitudes.
The estimated voltage sensitivity matrices, $\bm{R}$ and $\bm{X}$, are computed according to \eqref{eq:RX}.
The topology can be estimated using a variety of methods (see, e.g., \cite{Topology1, Topology2, Topology3}).
In this paper, we assume the node-to-edge incidence matrix is known, then only the line impedances need to be estimated.
After that, the estimated $\bm{R}$ and $\bm{X}$, denoted respectively by $\hat{\bm{R}}$ and $\hat{\bm{X}}$, are sent to the voltage controller.
The voltage controller then computes the set-points for the DER power injections that minimize some cost function subject to constraints \textbf{C1}--\textbf{C3}.
The DERs will be instructed to inject active and reactive power by the amount determined by the voltage controller.
A new set of measurements will be available once the DERs have modified their power injections.
These measurements will be used by the estimator to update $\hat{\bm{R}}$ and $\hat{\bm{X}}$ so as to reflect any changes in them.
The detailed formulations for the voltage sensitivity problem and the voltage control problem are presented next.

\subsection{Voltage Sensitivity Estimator} \label{sec:vse}

The voltage sensitivities can be estimated based on the LinDistFlow model in \eqref{eq:LinDistFlow}.
Assume at time instant $K$, we have measurements $V_0^{(k)}$, $\bm{V}^{(k)}$, $\bm{p}^{(k)}$, $\bm{q}^{(k)}$, $k \in \mathcal{K} = \{1, 2, \cdots, K\}$, where the index $k$ indicates the corresponding measurement is obtained at time instant $k$.
To reduce the computational burden, we select a subset of measurements from $\mathcal{K}$, denoted by $\mathcal{K}^\star = \{k_1, \cdots, k_m\}$.
The objective of the voltage sensitivity estimator at time instant $K$ is to estimate $\bm{R}$ and $\bm{X}$---essentially $\bm{r}$ and $\bm{x}$---from the selected measurements in $\mathcal{K}^\star$.

By using these measurements and the relation in \eqref{eq:LinDistFlow}, we can formulate the voltage sensitivity estimation problem as
\begin{equation} \label{eq:vse}
	\minimize_{\bm{r}, \bm{x}} \sum_{k \in \mathcal{K}^\star} \norm{\bm{R} \bm{p}^{(k)} + \bm{X} \bm{q}^{(k)} - \tilde{\bm{u}}^{(k)}}_2^2,
\end{equation}
where $\norm{\cdot}_2$ denotes the Euclidean norm.
We next show that \eqref{eq:vse} has a closed-form solution.

The matrix $\diag{\bm{x}}$ can be represented as follows:
\begin{equation} \label{eq:basis_expan}
	\diag{\bm{x}} = \sum_{i=1}^L x_i \bm{e_i} \bm{e_i}\trans,
\end{equation}
where $\bm{e}_i$ is the $i^{\mathrm{th}}$ basis vector in $\real^L$, i.e., all entries in $\bm{e}_i$ are $0$ except the $i^{\mathrm{th}}$ entry, which equals to $1$.
Using \eqref{eq:basis_expan}, we obtain that
\begin{equation}
\begin{array}{lll}
	\bm{X} \bm{q}^{(k)} & = & 2 (\bm{M}^{-1})\trans \diag{\bm{x}} \bm{M}^{-1} \bm{q}^{(k)} \\
	& = & 2 (\bm{M}^{-1})\trans \sum \limits_{i=1}^L x_i \bm{e_i} \bm{e_i}\trans \bm{M}^{-1} \bm{q}^{(k)} \\
	& = & \sum \limits_{i=1}^L \bm{\Gamma}_i \bm{q}^{(k)} x_i,
\end{array}	
\end{equation}
where $\bm{\Gamma}_i = 2 (\bm{M}^{-1})\trans \bm{e_i} \bm{e_i}\trans \bm{M}^{-1}$.
Similarly,
\begin{equation}
	\bm{R} \bm{p}^{(k)} = \sum \limits_{i=1}^L \bm{\Gamma}_i \bm{p}^{(k)} r_i = \sum \limits_{i=1}^L \bm{\Gamma}_i \alpha_i \bm{p}^{(k)} x_i.
\end{equation}

Define $\bm{\phi} = [(\tilde{\bm{u}}^{(k_1)})\trans, \cdots, (\tilde{\bm{u}}^{(k_m)})\trans]\trans$, and
\begin{equation*}
	\bm{\Phi} = 
	\begin{bmatrix}
		\bm{\Gamma}_1 (\alpha_1 \bm{p}^{(k_1)}+\bm{q}^{(k_1)}) & \cdots & \bm{\Gamma}_L (\alpha_L \bm{p}^{(k_1)}+\bm{q}^{(k_1)}) \\
		\vdots & \vdots & \vdots \\
		\bm{\Gamma}_1 (\alpha_1 \bm{p}^{(k_m)}+\bm{q}^{(k_m)}) & \cdots & \bm{\Gamma}_L (\alpha_L \bm{p}^{(k_m)}+\bm{q}^{(k_m)}) 
	\end{bmatrix}.
\end{equation*}
Note that both $\bm{\Phi} \in \real^{mL \times L}$ and $\bm{\phi} \in \real^{mN}$ are dependent on $\mathcal{K}^\star$.
Then \eqref{eq:vse} can be equivalently formulated in the classical form of a linear regression problem as follows:
\begin{equation} \label{eq:vse2}
	\minimize_{\bm{x}} \norm{\bm{\Phi} \bm{x} - \bm{\phi}}_2^2.
\end{equation}
When $\rank{\bm{\Phi}} = N$, the solution to \eqref{eq:vse2} is given by
\begin{equation} \label{eq:vse_sol}
	\hat{\bm{x}} = (\bm{\Phi}\trans \bm{\Phi})^{-1} \bm{\Phi}\trans \bm{\phi}.
\end{equation}
The resistance vector can be estimated using $\hat{\bm{r}} = \diag{\bm{\alpha}} \hat{\bm{x}}$.
Then, estimates of the voltage sensitivity matrices, denoted by $\hat{\bm{R}}$ and $\hat{\bm{X}}$, are computed according to \eqref{eq:RX}. 
However, due to the collinearity in the measurements \cite{ZJB}, $\bm{\Phi}\trans \bm{\Phi}$ may be ill-conditioned.
To obtain a numerically more stable solution to \eqref{eq:vse2}, the singular value decomposition technique is utilized to get the pseudo-inverse of $\bm{\Phi}$, denoted by $\bm{\Phi}^\dagger$; then the solution to \eqref{eq:vse2} is $\hat{\bm{x}} =  \bm{\Phi}^\dagger \bm{\phi}$.
Once new measurements are available at time instant $K+1$, the oldest measurements in $\mathcal{K}^\star$, i.e., the ones indexed by $k_1$ are replaced by the new one.

\subsection{Voltage Controller} \label{sec:vc}

The voltage controller aims to determine the optimal set-points for the DER power injections while meeting all requirements discussed in Section \ref{sec:vcp}.
Note that for a given set of power injections, the resulting voltage magnitude at each bus can be estimated using \eqref{eq:LinDistFlow}, where $\hat{\bm{R}}$ and $\hat{\bm{X}}$ are used instead of $\bm{R}$ and $\bm{X}$.
Denote the objective function as $c(\bm{p}^g, \bm{q}^g, \bm{u})$; then the voltage control problem can be formulated as the following convex program:
\begin{equation*} 
\minimize_{\bm{p}^g, \bm{q}^g} ~ c(\bm{p}^g, \bm{q}^g, \bm{u})
\end{equation*}
subject to 
\begin{subequations} \label{eq:vcp}
\begin{align} \label{eq:vcp_c1}
    \bm{u} - u_0 \ones_N = \hat{\bm{R}} (\bm{p}^g - \bm{p}^d) + \hat{\bm{X}} (\bm{q}^g - \bm{q}^d),
\end{align}
\begin{equation} \label{eq:vcp_c2}
	\underline{\bm{p}}^g \leq \bm{p}^g \leq \overline{\bm{p}}^g,
\end{equation}
\begin{equation} \label{eq:vcp_c3}
	\underline{\bm{q}}^g \leq \bm{q}^g \leq \overline{\bm{q}}^g,
\end{equation}
\begin{equation} \label{eq:vcp_c4}
	-\overline{\bm{f}} \leq \bm{M}^{-1} \bm{p} \leq \overline{\bm{f}},
\end{equation}
\end{subequations}
with
\begin{equation} \label{eq:obj}
	c(\bm{p}^g, \bm{q}^g, \bm{u}) = \bm{1}_N\trans (\bm{p}^g + \bm{q}^g) + \gamma ([\underline{\bm{u}} - \bm{u}]_+ + [\bm{u} - \overline{\bm{u}}]_+),
\end{equation}
where $\underline{\bm{u}} = [\underline{V}_1^2, \cdots, \underline{V}_N^2]\trans$, $\overline{\bm{u}} = [\overline{V}_1^2, \cdots, \overline{V}_N^2]\trans$, $[~]_+$ returns a non-negative vector, and $\gamma \in \real_+$. 
The first term of $c(\cdot)$ is the sum of active and power injections, while the second term penalizes the violation of constraint \textbf{C2}.

Constraint \eqref{eq:vcp_c1} is the LinDistFlow model, which is used to predict the voltage magnitudes for given power injections.
Note that $\bm{p}^d$ and $\bm{q}^d$ are measured before solving the voltage control problem.
Constraints \eqref{eq:vcp_c2} and \eqref{eq:vcp_c3} are essentially constraint \textbf{C1}.
Constraint \eqref{eq:vcp_c4} is essentially constraint \textbf{C3}.
Constraint \textbf{C2} is relaxed and included in the $c(\cdot)$.
Solving \eqref{eq:vcp} gives the optimal set-points for the DER power injections.

\section{Numerical Simulations} \label{sec:sim}

\begin{figure}[!t]
\centering
\includegraphics[width=3in]{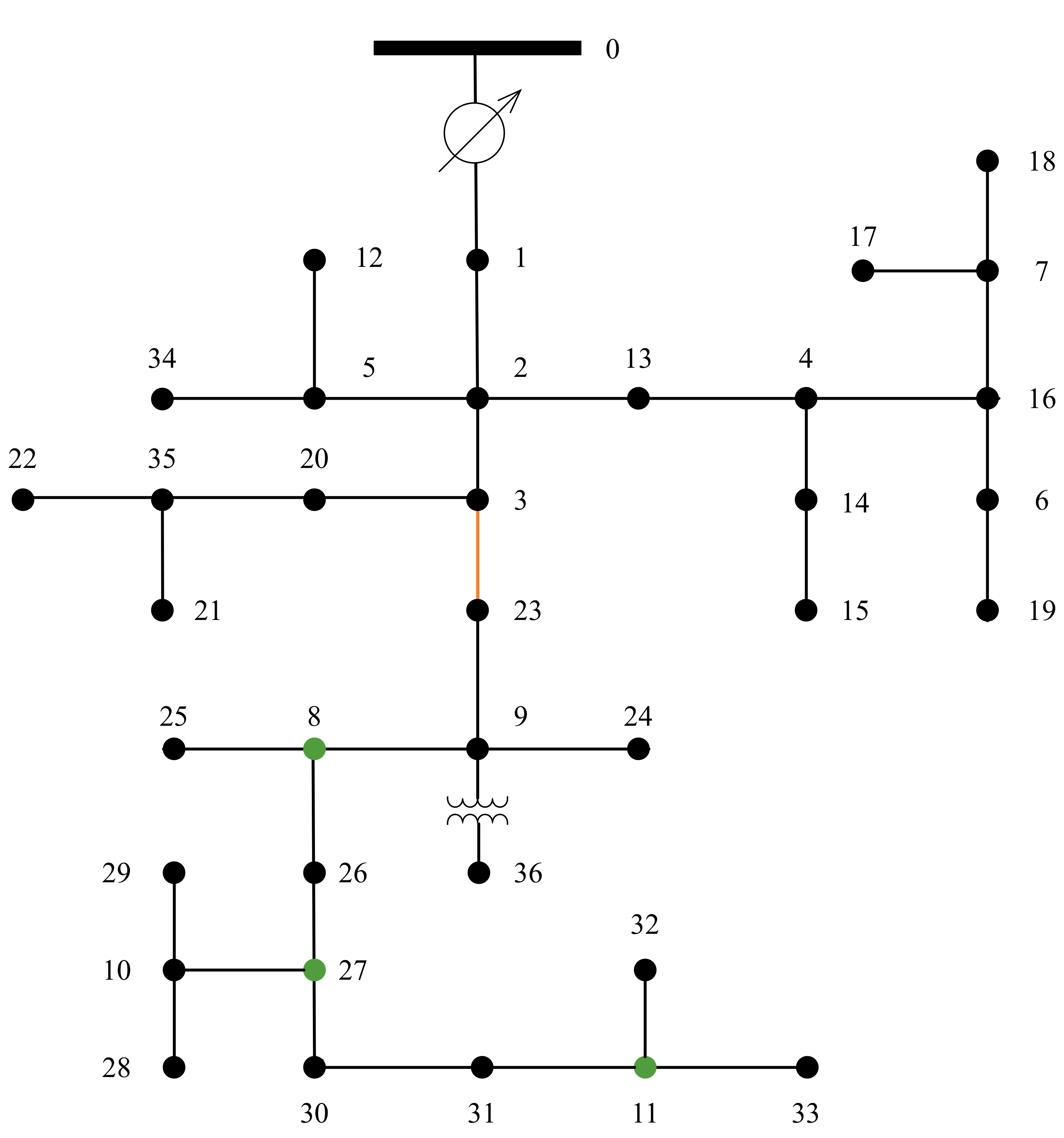}
\caption{IEEE 37-bus distribution test feeder.}
\label{fig:37bus}
\end{figure}

In this section, we illustrate the application of the proposed voltage control framework.
The parameter estimation accuracy and the voltage regulation performance are investigated.
A modified single-phase IEEE 37-bus distribution test feeder from \cite{test_feeder}, the topology of which is shown in Fig. \ref{fig:37bus}, is used for all numerical simulations.
Due to space limitation, the impacts of line flow constraints are not presented.
DERs are added at buses 8, 11, 27, with no active power, and maximum reactive power of 200, 300, and 250~kVAr, respectively.
The active load at bus $i$ is simulated by $p_i^d = p_i^{d0} (1 + \sigma^d \nu_i^{(k)})$, where $p_i^{d0}$ is the nominal active load at bus $i$, $\sigma^d = 0.01$, $\nu_i^{(k)}$ is a standard Gaussian random variable.
The reactive load is simulated in a similar manner.
A white noise with standard deviation of $2\times 10^{-4}$ is also added as measurement error \cite{VS1}.
The minimum and maximum voltages are 0.95 p.u. and 1.05 p.u., respectively.

\subsection{Estimation Accuracy}

\begin{figure}[!t]
\centering
\includegraphics[width=3.5in]{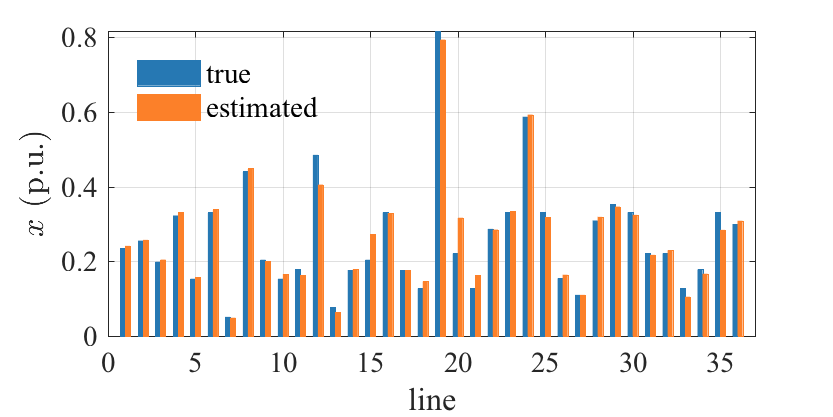}
\caption{Line parameter estimation results.}
\label{fig:line_params_x}
\end{figure}

The estimation accuracy is evaluated using the mean absolute error (MAE).
Figure \ref{fig:line_params_x} shows the estimated line parameters when 20 measurements are utilized; the MAE of $\bm{x}$ is $0.0160$~p.u., and that of $\bm{X}$ is $0.0406$~p.u, both of which are relatively small compared to the corresponding mean values of $\bm{x}$ and $\bm{X}$, which are $0.2622$~p.u. and $1.6812$~p.u., respectively.
The results for $\bm{r}$ and $\bm{R}$ are similar and are thus not presented due to space limitation.

\begin{figure}[!t]
\centering
\includegraphics[width=3.5in]{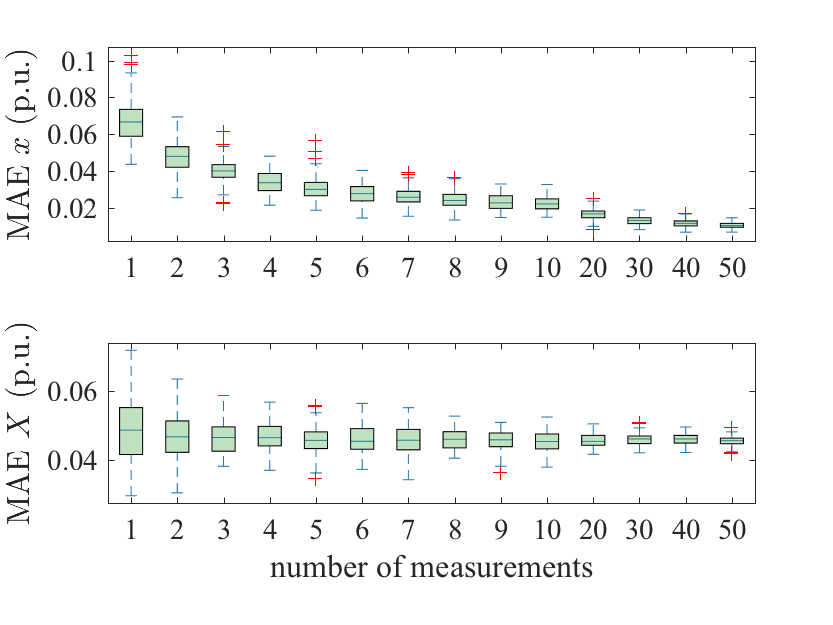}
\caption{Line parameter estimation accuracy versus number of measurements.}
\label{fig:line_param_error}
\end{figure}

Theoretically, line parameters can be estimated using as few as one set of measurements because the network is radial; thus, $\bm{\Phi}$ in \eqref{eq:vse2} would be $(L\times L)$-dimensional matrix.
It should be noted that increasing the number of measurements does improve the estimation accuracy.
Figure \ref{fig:line_param_error} presents a box plot of estimation accuracy with respect to the number of measurements.
For a given number of measurements, the line parameters are estimated in 100 Monte Carlo simulation runs.
Both the mean and variance of the MAE of $\bm{x}$ decrease when the number of measurements increases. 
While the mean of the MAE of $\bm{X}$ does not change much, its variance decreases significantly.

\subsection{Voltage Control Effects}

\begin{figure}[!t]
\centering
\includegraphics[width=3.5in]{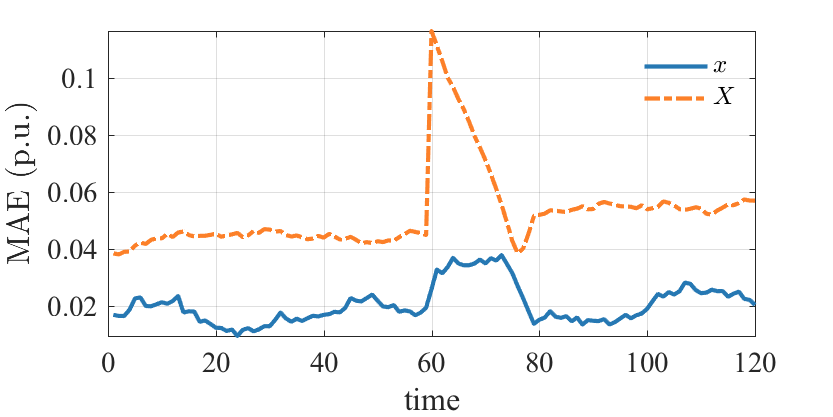}
\caption{Evolution of MAEs.}
\label{fig:MAE2}
\end{figure}

\begin{figure}[!t]
\centering
\includegraphics[width=3.5in]{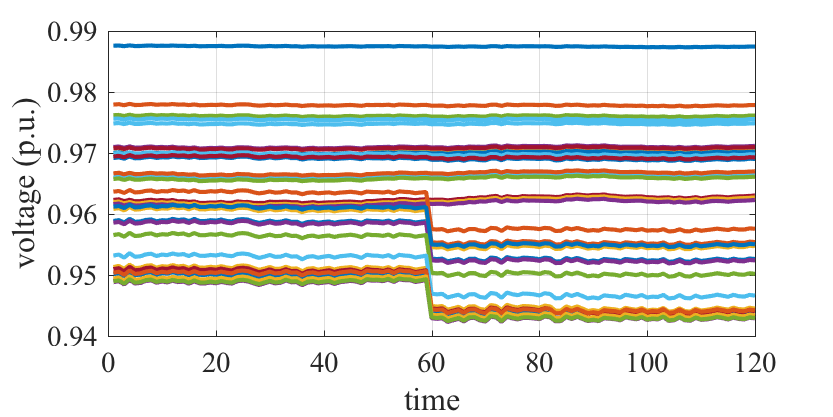}
\caption{Voltage profiles under model-based control.}
\label{fig:voltage1}
\end{figure}

\begin{figure}[!t]
\centering
\includegraphics[width=3.5in]{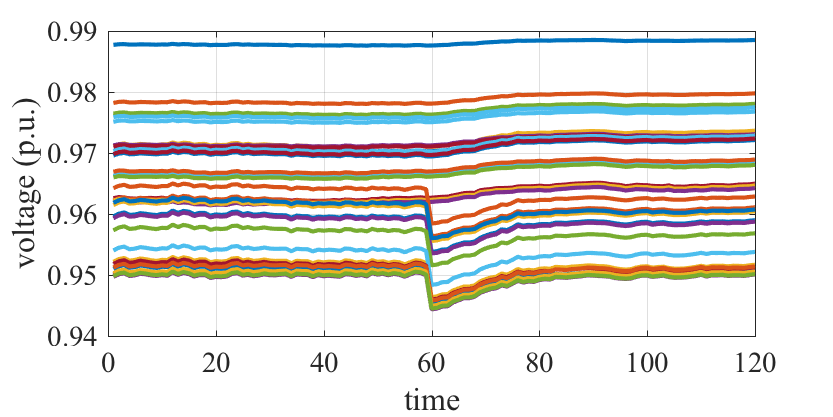}
\caption{Voltage profiles under proposed control .}
\label{fig:voltage2}
\end{figure}

\begin{figure}[!t]
\centering
\includegraphics[width=3.5in]{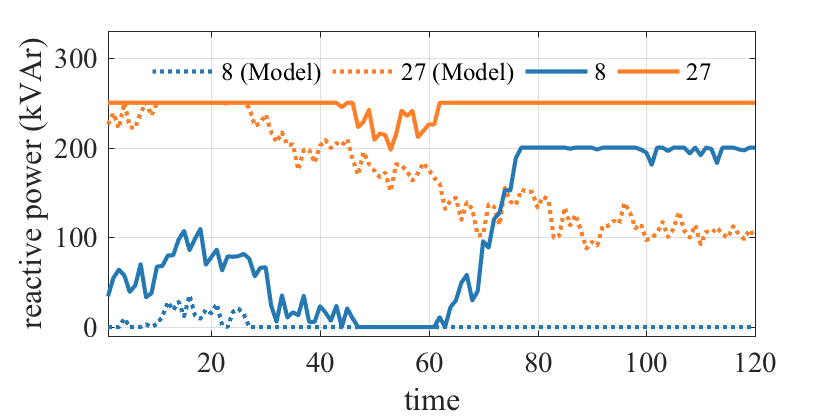}
\caption{Evolution of reactive power profiles.}
\label{fig:reactive_power_comp}
\end{figure}

Data-driven methods are inherently adaptive to changes in system conditions given that sufficient amount of data are obtained timely \cite{MB1, MB2}.
This, however, may result in high communication burden. 
A distinct advantage of the proposed estimator is that it requires far fewer data compared to other sensitivity estimators by exploiting the structural characteristics, allowing it to track the changes faster.
Consequently, the voltage controller, which is built upon the estimator, can adapt to system changes quickly.

To illustrate the adaptivity of the propose framework, we perturb the line parameters at line (3,~23) at time instant 60.
In the simulation, the number of measurements used by the estimator is 20. 
The value of the parameter $\gamma$ in \eqref{eq:obj} is set to 1000.
As shown in Fig. \ref{fig:MAE2}, after the perturbation, the MAEs of $\bm{x}$ and $\bm{X}$ increased, but were reduced soon.
As shown in Fig. \ref{fig:voltage1}, when the true model parameters before perturbation are used in the voltage controller, the voltages dropped immediately after the perturbation but did not return to the pre-specified range.
In contrast, when the estimated parameters are used, the voltages were restored to the pre-specified range after 20 time instances as shown in Fig. \ref{fig:voltage2}.
Figure \ref{fig:reactive_power_comp} shows the DER reactive power outputs; note that the reactive power of DER 11 is 300~kVAr throughout the simulated time interval.

\section{Concluding Remarks} \label{sec:con}

In this paper, we proposed an adaptive data-driven voltage control framework for power distribution systems, which consists of a voltage sensitivity estimator and a voltage controller. In particular, the proposed estimator requires much fewer data than existing ones by exploiting the structural characteristics of the power distribution system. The inherent nature of the estimator makes it adaptive to changes in system conditions. Numerical simulation shows that the framework is effective and efficient in coordinating the DERs to provide voltage control. Future work will develop an efficient topology estimator using the same set of data so as to complete the framework.

\bibliographystyle{IEEEtran}
\bibliography{ref}

\begin{thebibliography}{10}
\providecommand{\url}[1]{#1}
\csname url@samestyle\endcsname
\providecommand{\newblock}{\relax}
\providecommand{\bibinfo}[2]{#2}
\providecommand{\BIBentrySTDinterwordspacing}{\spaceskip=0pt\relax}
\providecommand{\BIBentryALTinterwordstretchfactor}{4}
\providecommand{\BIBentryALTinterwordspacing}{\spaceskip=\fontdimen2\font plus
\BIBentryALTinterwordstretchfactor\fontdimen3\font minus
  \fontdimen4\font\relax}
\providecommand{\BIBforeignlanguage}[2]{{%
\expandafter\ifx\csname l@#1\endcsname\relax
\typeout{** WARNING: IEEEtran.bst: No hyphenation pattern has been}%
\typeout{** loaded for the language `#1'. Using the pattern for}%
\typeout{** the default language instead.}%
\else
\language=\csname l@#1\endcsname
\fi
#2}}
\providecommand{\BIBdecl}{\relax}
\BIBdecl

\bibitem{VC1}
B.~A. Robbins, C.~N. Hadjicostis, and A.~D. Dom\'inguez-Garc\'ia, ``A two-stage
  distributed architecture for voltage control in power distribution systems,''
  \emph{IEEE Trans. Power Syst.}, vol.~28, no.~2, pp. 1470--1482, May 2013.

\bibitem{VC2}
B.~A. Robbins and A.~D. Dom\'inguez-Garc\'ia, ``Optimal reactive power dispatch
  for voltage regulation in unbalanced distribution systems,'' \emph{IEEE
  Trans. Power Syst.}, vol.~31, no.~4, pp. 2903--2913, July 2016.

\bibitem{MB1}
Y.~C. Chen, A.~D. Dom\'inguez-Garc\'ia, and P.~W. Sauer, ``Measurement-based
  estimation of linear sensitivity distribution factors and applications,''
  \emph{IEEE Trans. Power Syst.}, vol.~29, no.~3, pp. 1372--1382, May 2014.

\bibitem{MB2}
Y.~C. Chen, S.~V. Dhople, A.~D. Dom\'inguez-Garc\'ia, and P.~W. Sauer,
  ``Generalized injection shift factors,'' \emph{IEEE Trans. on Smart Grid},
  vol.~8, no.~5, pp. 2071--2080, Sept 2017.

\bibitem{MB3}
Y.~C. Chen, A.~D. Dom\'inguez-Garc\'ia, and P.~W. Sauer, ``A sparse
  representation approach to online estimation of power system distribution
  factors,'' \emph{IEEE Trans. Power Syst.}, vol.~30, no.~4, pp. 1727--1738,
  July 2015.

\bibitem{MB4}
K.~E.~V. Horn, A.~D. Dom\'inguez-Garc\'ia, and P.~W. Sauer, ``Measurement-based
  real-time security-constrained economic dispatch,'' \emph{IEEE Trans. Power
  Syst.}, vol.~31, no.~5, pp. 3548--3560, Sept 2016.

\bibitem{ZJB}
J.~Zhang, X.~Zheng, Z.~Wang, L.~Guan, and C.~Y. Chung, ``Power system
  sensitivity identification: Inherent system properties and data quality,''
  \emph{IEEE Trans. Power Syst.}, vol.~32, no.~4, pp. 2756--2766, July 2017.

\bibitem{33bus}
M.~E. Baran and F.~F. Wu, ``Network reconfiguration in distribution systems for
  loss reduction and load balancing,'' \emph{IEEE Trans. Power Del.}, vol.~4,
  no.~2, pp. 1401--1407, Apr 1989.

\bibitem{DistFlow}
V.~Kekatos, L.~Zhang, G.~B. Giannakis, and R.~Baldick, ``Voltage regulation
  algorithms for multiphase power distribution grids,'' \emph{IEEE Trans. Power
  Syst.}, vol.~31, no.~5, pp. 3913--3923, Sept 2016.

\bibitem{Topology1}
S.~Bolognani, N.~Bof, D.~Michelotti, R.~Muraro, and L.~Schenato,
  ``Identification of power distribution network topology via voltage
  correlation analysis,'' in \emph{Proc. of IEEE Conference on Decision and
  Control}, Dec 2013, pp. 1659--1664.

\bibitem{Topology2}
Y.~Weng, Y.~Liao, and R.~Rajagopal, ``Distributed energy resources topology
  identification via graphical modeling,'' \emph{IEEE Trans. Power Syst.},
  vol.~32, no.~4, pp. 2682--2694, July 2017.

\bibitem{Topology3}
D.~Deka, S.~Backhaus, and M.~Chertkov, ``Estimating distribution grid
  topologies: A graphical learning based approach,'' in \emph{Proc. of Power
  Systems Computation Conference}, June 2016, pp. 1--7.

\bibitem{test_feeder}
\BIBentryALTinterwordspacing
{IEEE} distribution test feeders. [Online]. Available:
  \url{https://ewh.ieee.org/soc/pes/dsacom/testfeeders/}
\BIBentrySTDinterwordspacing

\bibitem{VS1}
C.~Mugnier, K.~Christakou, J.~Jaton, M.~D. Vivo, M.~Carpita, and M.~Paolone,
  ``Model-less/measurement-based computation of voltage sensitivities in
  unbalanced electrical distribution networks,'' in \emph{Proc. of Power
  Systems Computation Conference}, June 2016, pp. 1--7.

\end{thebibliography}

\end{document}